\documentclass[11pt]{amsart}
\usepackage{amscd,amssymb}
\newcommand{\bdes}{\begin{description}}
\newcommand{\edes}{\end{description}}
\newcommand{\beqn}{\begin{eqnarray*}}
\newcommand{\eeqn}{\end{eqnarray*}}
\newcommand{\EE }{{\mathbb E}}
\newcommand{\PP }{{\mathbb P}}

\newtheorem{proposition}{Proposition}[subsection]
\newtheorem{lemma}{Lemma}[subsection]

\newtheorem{def/th}{Definition/Theorem}[subsection]

\begin{document}
\title{The Aspinwall-Morrison calculation and Gromov-Witten theory}
\author{Artur Elezi}
\date{\today}
\begin{abstract} We connect the Aspinwall-Morrison calculation to Gromov-Witten theory.
\end{abstract}
\maketitle
\section{\bf Introduction}
{I. A bit of history.} (See \cite{[5]} for a good reference on the history of the problem.) One of the problems in the old and recent story of mirror symetry has been the issue of multiple covers on a Calabi-Yau 3-fold $X$. In the pre Gromov-Witten era, this problem can be explained in terms of topological field theories.

Let $X$ be a Calabi-Yau threefold and $H_1,H_2,H_3\in H^2(X)$. The corresponding 3-point correlator in the A-model of $X$ is a path integral that can be expressed as follows:
\begin{equation}
\langle H_1,H_2,H_3 \rangle=\int_{X}H_1H_2H_3+\sum_{\beta\in H_2(X)}N_{\beta}(H_1,H_2,H_3)q^{\beta}.
\end{equation}
We explain the notation. The parameter $q=(q_1,...q_k)$ is a local coordinate on the K$\ddot{a}$hler moduli space of $X$. Let $(d_1,...,d_k)$ be the coordinates of $\beta$ with respect to an integral base of the Mori cone of $X$. Then $q^{\beta}:=q_1^{d_1}\cdot \cdot \cdot q_k^{d_k}$. 

The path integral is not a well defined notion, but beyond that, and probably more importantly, there is no rigorous definition of $N_{\beta}(H_1,H_2,H_3)$ in the framework of topological field theories. Let $Z_i$ for $i=1,2,3$ be a cycle whose fundamental class is Poincar$\acute{e}$ dual to $H_i$. Heuristically, the ``invariant'' $N_{\beta}(H_1,H_2,H_3)$ is described as the ``number'' of holomorphic maps in 
\begin{equation}
\{f:\PP^1\rightarrow X~|~f_*([\PP^1])=\beta,f(0)\in Z_1, f(1)\in Z_2, f(\infty)\in Z_3\}. 
\end{equation}
This is certainly not precise, for there may be infinitely many such maps. Let $C\subset X$ be a smooth rational curve. Fix an isomorphism $g:\PP^1\rightarrow C$. For any degree $k$ multiple cover $f:\PP^1\rightarrow \PP^1$ the composition $f\circ g:\PP^1\rightarrow C$ satisfies $(f\circ g)_*([\PP^1])=k[C]$. One would then naturally ask:

What is the contribution of the space of degree $k$ multiple covers of $C$ to the ``invariant'' $N_{k[C]}(H_1,H_2,H_3)$?

Since this question is about the numbers $N_{k[C]}(H_1,H_2,H_3)$, it is not a well posed one. It can be made precise in the framework of Gromov-Witten theory.

The answer was conjectured in \cite {[4]} by looking at the classical example of a Calabi-Yau. If $X$ is  a quintic threefold then $H^2(X)$ is one dimensional. Let $H$ be its generator. The 3-point correlator of the quintic can be calculated explicitly:
\begin{equation}
\langle H,H,H \rangle=5+\sum_{d=1}^{\infty}n_dd^3\frac{q^d}{1-q^d},
\end{equation}
where $n_d$ is the virtual number of degree $d$ rational curves (instantons) in the quintic. The instanton number $n_d$ agrees with the number of degree $d$ rational curves in the quintic if every rational curve of degree $d$ is smooth, isolated and with normal bundle $N=\mathcal O(-1)\oplus \mathcal O(-1)$. This is not the case for there are $6$-nodal rational plane quintic curves on a generic quintic threefold (see \cite {[12]}), hence a rigorous definition of the instanton numbers $n_d$ did not exist.

The last equation can be transformed as follows:
\begin{equation}
\langle H,H,H \rangle=5+\sum_{d=1}^{\infty}(\sum_{k|d}n_kk^3)q^d.
\end{equation}
By comparing it to the equation (1) we can see that:
\begin{equation}
N_d(H,H,H)=\sum_{k|d}n_kk^3.
\end{equation}
It looks that each degree $k$ rational curve $C$ in the quintic 3-fold $X$ contributes by: 
\begin{equation}
\int_{C}H\cdot \int_{C}H\cdot \int_{C}H
\end{equation}
to $N_d(H,H,H)$ for any $d$ such that $k|d$. 

For a general Calabi-Yau $X$, the (pre Gromov-Witten) multiple cover formula can be formulated as follows:

Let $C\subset X$ be a smooth, rational curve such that $N_{C/X}=\mathcal O_C(-1)\oplus \mathcal O_C(-1)$. The contribution of degree $k$ multiple covers of $C$ in $N_{k[C]}(H_1,H_2,H_3)$ is: 
\begin{equation}
\int_{C}H_1\cdot \int_{C}H_2\cdot \int_{C}H_3.
\end{equation}

It was in this form that the multiple cover formula was proven by Aspinwall and Morrison in \cite{[1]} and by Voisin in \cite{[10]}.

It follows from the above equation that:
\begin{eqnarray}
& & N_{\beta}(H_1,H_2,H_3)=\sum_{\beta=d\gamma}n_{\gamma}\int_{\gamma}H_1 \int_{\gamma}H_2 \int_{\gamma}H_3 \nonumber \\ & & =(\sum_{\beta=d\gamma}n_{\gamma}d^{-3})\int_{\beta}H_1 \int_{\beta}H_2 \int_{\beta}H_3
\end{eqnarray}
where $n_{\gamma}$ is the virtual number (instantons) of rational curves of type $\gamma$ in $X$.

A rigorous definition of $N_{\beta}$ and $n_{\beta}$ requires a new conceptual framework which is now known as Gromov-Witten theory. Let $X$ be a smooth, projective manifold and $\beta\in H_2(X)$. There is a moduli stack $\overline M_{0,n}(X,\beta)$ which parametrizes pointed, stable maps of degree $\beta$. Universal properties of these maduli stacks imply the existence of several features:
\begin{eqnarray}
& & e=(e_1,e_2,...,e_n):\overline M_{0,n}(X,\beta)\rightarrow X^n,~ \pi_n:\overline M_{0,n}(X,\beta)\rightarrow \overline M_{0,n-1}(X,\beta) \nonumber \\ & & \pi:\overline M_{0,n}(X,\beta)\rightarrow \overline M_{0,0}(X,\beta),~\hat{\pi}:\overline M_{0,n}(X,\beta)\rightarrow \overline M_{0,n}.
\end{eqnarray}
The map $e$ evaluates the pointed, stable map at the marked points, $\pi_n$ forgets the last marked point and collapses the unstable components of the source curve, $\pi$ forgets the marked points and $\hat{\pi}$ forgets the map and stabilizes the pointed source curve. The expected dimension of $\overline M_{0,n}(X,\beta)$ is $\text{dim}~X+\int_{\beta}(-K_X)+n-3$. The dimension of the moduli stack of stable maps may be greater than the expected dimension. In this case, a Chow class of the expected dimension has been constructed. It plays the role of the fundamental class, hence it is called the virtual fundamental class and denoted by $[\overline M_{0,n}(X,\beta)]^{\text{vir}}$ (see \cite{[11]},{[13]}).

Let $X$ be a Calabi-Yau threefold and $H_1,H_2,H_3\in H^2(X)$. In the Gromov-Witten setting:
\begin{equation}
N_{\beta}(H_1,H_2,H_3):=\int_{[\overline M_{0,3}(X,\beta)]^{\text{virt}}}e_1^*(H_1)e_2^*(H_2)e_3^*(H_3).
\end{equation} 
The expected dimension of $\overline M_{0,0}(X,\beta)$ is zero. Let:
\begin{equation}
N_{\beta}:=\text{deg}([\overline M_{0,0}(X,\beta)]^{\text{virt}})
\end{equation}
By the divisor axiom:
\begin{equation}
N_{\beta}(H_1,H_2,H_3)=N_{\beta}\int_{\beta}H_1\int_{\beta}H_2\int_{\beta}H_3.
\end{equation}

Let $C\subset X$ be a smooth rational curve with $N_{C/X}=\mathcal O_C(-1)\oplus \mathcal O_C(-1)$. The moduli space $\overline M_{0,0}(X,d[C])$ contains a component of positive dimension, namely $\overline M_{0,0}(C,d)$. The dimension of this component is $2d-2$. Consider the following diagram:
\[ \begin{CD}
\overline M_{0,1}(C,d)@>e_1>>C \\
@VV \pi V \\
\overline M_{0,0}(C,d)
\end{CD}\]

The sheaf:
\begin{equation}
V_d:=R^1\pi_*(\mathcal O_C(-1)\oplus \mathcal O_C(-1))
\end{equation}
is locally free of rank $2d-2$. Let $\EE_d$ be its top chern class. An assertion of Kontsevich in \cite{[6]}, which was proven by Behrend in \cite{[2]}, states that the part of $[\overline M_{0,0}(X,\beta)]^{\text{virt}}$ supported in $\overline M_{0,0}(C,d)$ is Poincar$\acute{e}$ dual to $\EE_d$. The multiple cover formula in this context says that: 
\begin{equation}
\int_{\overline M_{0,0}(C,d)}\EE_d=\displaystyle d^{-3}
\end{equation} 
i.e. the curve $C$ contributes by $d^{-3}$ to to $N_{d[C]}$. 

The multiple cover formula in this form was proven by Kontsevich \cite{[6]}, Lian-Liu-Yau \cite{[7]}, Manin \cite{[8]} and Pandharipande \cite{[9]}. 

By the divisor property, the multiple cover formula in this context follows from:
\begin{equation}
\int_{\overline M_{0,3}(C,d)}e_1^*(h)e_2^*(h)e_3^*(h)\pi^*(\EE_d)=1
\end{equation}
The instanton numbers $n_{\gamma}$ are defined inductively by:
\begin{equation}
N_{\beta}=\sum_{\beta=k\gamma}n_{\gamma}k^{-3}
\end{equation}
The point of this introduction is that the Aspinwall-Morrison calculation deals with concepts and questions that were not well defined at the time. Hence their calculation, although useful and convincing, is incomplete. The purpose of this paper is to relate the two calculations, hence justifying the Aspinwall-Morrison calculation and closing this historic chapter in the subject.

We show in passing the connection between the two formulations of the multiple cover formula for the quintic threefold:
\begin{equation}
N_d(H,H,H)=d^3N_d=d^3\sum_{k|d}n_k(\frac{k}{d})^3=\sum_{k|d}n_kk^3
\end{equation}

{\bf II. A review of the Aspinwall-Morrison calculation.} Consider a Calabi-Yau threefold $X$ and a rational curve $C\subset X$ such that $N_{C/X}=\mathcal O_C(-1)\oplus \mathcal O_C(-1)$. Let: 
\begin{equation}
N_d(C):=\{f:\PP^1\rightarrow X ~|~ f(\PP^1)=C, \text{deg}f=d\}
\end{equation}
be the space of parameterized maps from $\PP^1$ to $X$. Since $C$ is isolated, $N_d(C)$ is a component of the space of all maps from $\PP^1$ to $X$. 

At a moduli point $[f]$, the tangent space and the obstruction space are given respectively by $H^0(f^*(T_X))$ and $H^1(f^*(T_X))$, i.e. locally $N_d(C)$ is given by dim $H^1(f^*(T_X))$ equations in the tangent space. The virtual dimension is:
\begin{equation}
\text{dim}~H^0(f^*(T_X))-\text{dim}~H^1(f^*(T_X))=3.
\end{equation}

The space $N_d(C)$ compactifies to $\overline N_d(C)=\PP^{2d+1}$. Let $\overline \Gamma$ be the compactification of the universal graph $\Gamma\subset N_d(C)\times \PP^1\times C$ and $H$ the hyperplane class in $\overline N_d(C)$.. 

The dimension of $H^1(f^*(T_X))$ is $2d-2$ for any $f$. These vector spaces fit together to form a bundle $\mathcal U_d$ over $N_d(C)$. Let $p_i$ be the $i$-th projection on $\overline N_d(C)\times \PP^1\times C$. The bundle $\mathcal U_d$ extends to: 
\begin{equation}
U_d:=R^1{p_1}_*(p_3^*(T_X|C)|_{\overline \Gamma})
\end{equation}
over $\overline N_d(C)$. A calculation in \cite{[1]} yields $U_d=\mathcal O(-1)^{\oplus {d-2}}$. Based primarily on considerations from topological field theories, Aspinwall and Morrison asserted that the cycle corresponding to the degree $d$ multiple covers of $C$ is Poincar$\acute{e}$ dual to $c_{\text{top}}(U_d)=H^{2d-2}$. We will see that this is consistent with the notion of the virtual fundamental class.

Let $H_i\in H^2(X)$ for $i=1,2,3$ and $Z_i$ their Poincar$\acute{e}$ duals. The space: 
\begin{equation}
\{f\in N_d(C)~|~f(0)=0\}
\end{equation}
gives rise to a linear subspace of $\overline N_d(C)$. Therefore:
\begin{eqnarray}
& & \#\{f\in N_d(C)~|~f(0)=0,f(1)=1,f(\infty)=\infty\}\nonumber \\ & & =\int_{\overline N_d(C)}H\cdot H\cdot H\cdot c_{\text{top}}U_d=1.
\end{eqnarray}

It follows that the contribution of $N_d(C)$ to: 
\begin{equation}
\#\{f:\PP^1\rightarrow X~|~f_*[\PP^1]=d[C],f(0)\in Z_1, f(1)\in Z_2, f(\infty)\in Z_3\}
\end{equation}
is 
\begin{equation}
\int_{C}H_1\cdot \int_{C}H_2\cdot \int_{C}H_3.
\end{equation}

We emphasize that the multiple cover formula in this approach follows from:
\begin{equation}
\int_{\overline N_d(C)}H\cdot H\cdot H\cdot c_{\text{top}}U_d=\int_{\overline N_d(C)}H^{2d+1}=1.
\end{equation}

{\bf III. The connection to the Gromov-Witten theory}. The main result in this paper is the following: 
\begin{proposition} There exists a birational morphism:
\begin{equation}
\alpha:\overline M_{0,3}(C,d)\rightarrow {\overline N}_d(C)
\end{equation}
such that:
\begin{enumerate}
\item $\alpha_*(e_i^*(h))=H$ for $i=1,2,3.$
\item $\alpha_*(e_1^*(h)e_2^*(h)e_3^*(h))=H^3$
\item $\alpha_*(e_1^*(h)e_2^*(h)e_3^*(h)\pi^*(\EE_d))=H^{2d+1}$.
\end{enumerate}
\end{proposition}
This proposition implies that the equations $(15)$ and $(25)$ are equivalent, hence connecting the Aspinwall-Morrison calculation to the Gromov-Witten theory.

{\bf Acknowledgements}. The problem was suggested to the author by Sheldon Katz (see also the note in \cite{[3]}) who was very helpful through this work. We would also like to thank Jun Li for fruitful discussions on the subject.
\section{\bf Relation of the Aspinwall-Morrison formula with Gromov-Witten invariants}
The space of nonparameterized degree $d$ maps $f:\PP^1\rightarrow \PP^n$ has two particular compactifications that have been employed successfully especially in proving mirror theorems for projective spaces: the nonlinear sigma model (or the graph space): 
\begin{equation}
M^n_d:=\overline M_{0,0}(\PP^n\times \PP^1,(d,1))
\end{equation}
and the linear sigma model: 
\begin{equation}
N^n_d:=\PP(H^0(\mathcal O_{\PP^1}(d))).
\end{equation} 
Elements of $N^n_d$ are $(n+1)$-tuples $[P_0,...,P_n]$ of degree $d$ polynomials in two variables $w_0,w_1$. The linear sigma model $N_d$ is a projective space via the identification $[P_0,...]=[\sum_{i}a_iw_0^iw_1^{d_i},...]=[a_0,...,a_d,...]$. Note that $N^1_d=\overline N_d(C)$ for $C\simeq \PP^1$. Let $H$ be the hyperplane class in $N^n_d$. 

There exists a birational morphism $\phi:M^n_d\rightarrow N^n_d$. We describe this morphism set-theoretically. Let $(C',f)\in M^n_d$. There is a unique component $C_0$ of $C'$ that is mapped with degree $1$ to $\PP^1$. Let $C_1,...,C_r$ be the irreducible components of the rest of the curve and $q_i=(c_i,d_i)$ the nodes of $C'$ on $C_0$. Let $d-i$ be the degree of the map $p_2\circ f:C'\rightarrow \PP^n$ on $C_i$ for $i=0,1,...,r$. Let $R(w_0,w_1)=\prod_{i=1}^{r}(c_iw_1-d_iw_0)^{d_i}$. If the restriction of the map $p_2\circ f$ is given by $[Q_0,...,Q_n]$ then:
\begin{equation}
\phi(C',f):=[RQ_0,...,RQ_n].
\end{equation}
A proof of the fact that $\phi$ is a morphism is given by J. Li in \cite{[7]}.

The first step in connecting the Aspinwall-Morrison calculation to Gromov-Witten invariants is showing that $M^n_d$ and $N^n_d$ are birational models for $\overline M_{0,3}(\PP^n,d)$.

\begin{lemma} There exists a birational map $\psi:\overline M_{0,3}(\PP^n,d)\rightarrow M^n_d$.
\end{lemma}

{\bf Proof}. Consider the following diagram:
\[ \begin{CD}
\overline M_{0,4}(\PP^n,d)@>(\hat{\pi},e_4)>>\overline M_{0,4}\times \PP^n \\
@VV \pi_4 V  \\
\overline M_{0,3}(\PP^n,d).
\end{CD} \]
Since $\overline M_{0,4}\simeq \PP^1$ and $e_4$ is stable in the fibers of $\pi_4$, the above diagram exhibits a stable family of maps of degree $(1,d)$ parametrized by $\overline M_{0,3}(\PP^n,d)$. Universal properties of $M^n_d$ yield a morphism: 
\begin{equation}
\psi:\overline M_{0,3}(\PP^n,d)\rightarrow M^n_d.
\end{equation}
The map $\psi$ is an isomorphism in the smooth locus, hence it is a birational map.$\dagger$

\vspace{0.2in}

Let $\pi_4:\overline M_{0,4}\rightarrow \overline M_{0,3}=\{pt\}$ be the map that forgets the last marked point and $\sigma_i$ be the section of the i-th marked point for $i=1,2,3$. Choose coordinates on $\overline M_{0,4}\simeq \PP^1$ such that the images of these three sections are respectively $0=[1,0],\infty=[0,1],1=[1,1]$. Let 
\begin{equation}
\alpha:=\phi\circ \psi:\overline M_{0,3}(\PP^n,d)\rightarrow N^n_d.
\end{equation}

\begin{proposition} Let $h$ be the hyperplane class of $\PP^n$. 
\begin{enumerate}
\item $\alpha_*(e_i^*(h))=H$ for $i=1,2,3$.
\item $\alpha_*(e_1^*(h)e_2^*(h)e_3^*(h))=H^3$
\end{enumerate}
\end{proposition}

{\bf Proof}. Let 
\begin{equation}
\nu_1:N_d---> \PP^n
\end{equation}
be a rational map defined by 
\begin{equation}
\nu_1([P_0,P_1,...,P_n])=[P_0(1,0),P_1(1,0),...,P_n(1,0)]. 
\end{equation}
This map is defined in the complement $U$ of a codimension $n+1$ linear subspace $P(W_1)$ of $N^n_d$. Clearly $\nu_1^*(h)=H$ on $U$. The preimage $D_{1,23}$ of $P(W_1)$ in $\overline M_{0,3}(\PP^n,d)$ is a sum of $d$ boundary divisors $D(\{x_1\},\{x_2,x_3\},d_1,d_2)$ with $d_1>0$ and $d_1+d_2=d$. The evaluation map $e_1$ over $U$ factors through the rational map $\nu_1$. It follows that 
\begin{equation}
e_1^*(h)=\alpha^*(H)+D_1,
\end{equation} 
where $D_1$\footnote{It can be shown that $D_1=-\sum_{d_1}d_1D(\{x_1\},\{x_2,x_3\},d_1,d-d_1)$ but this is not important in this paper.} is a divisor supported in $D_{1,23}$. Using the evaluations at $1$ and $\infty$ on $N^n_d$, we obtain: 
\begin{equation}
e_2^*(h)=\alpha^*(H)+D_2
\end{equation}
and 
\begin{equation}
e_3^*(h)=\alpha^*(H)+D_3, 
\end{equation}
where $D_2$ is a divisor  supported in $D_{2,13}$ and $D_3$ is supported in $D_{3,12}$. 

The $\psi$-image of $D(\{x_1\},\{x_2,x_3\},d_1,d_2)$ does not detect the movement of the marking $x_1$ along its incident component, hence it is a codimension $2$ cycle in $M^n_d$. It follows that $\psi_*(D_1)=0$. Similarly $\psi_*(D_2)=0$ and $\psi_*(D_3)=0$. Both $\psi$ and $\phi$ are birational hence by the projection formula:
\begin{equation}
\alpha_*(e_i^*(h))=H
\end{equation}
for $i=1,2,3$.

Let $D'\in D_{1,23},D''\in D_{2,13},D'''\in D_{3,12}$ be irreducible boundary divisors. The intersection of any two of them either is $0$ or its image is a codimension $4$ cycle in $M^n_d$. It follows that: 
\begin{equation}
\psi_*(D'D'')=\psi_*(D'D''')=\psi_*(D''D''')=0.
\end{equation} 
Notice also that: 
\begin{equation}
D'D''D'''=0.
\end{equation} 
The projection formula yields:
\begin{equation}
\psi_*(e_1^*(h)e_2^*(h)e_3^*(h)=\psi_*(\prod_{i}(\psi^*(\phi^*(H))+D_i))=\prod_{i}(\phi^*(H))=\phi^*(H^3).
\end{equation}
The lemma follows from the fact that $\phi$ is a birational map.$\dagger$

\vspace{0.2in}

Return now to the case $n=1$ of our interest.

Let $\rho:M^1_d\rightarrow \overline M_{0,0}(C,d)$ be the natural morphism. The composition: 
\begin{equation}
\rho\circ \psi:\overline M_{0,3}(C,d)\rightarrow \overline M_{0,0}(C,d)
\end{equation}
 is the map $\pi$ that forgets the $3$ marked points and stabilizes the source curve. Recall Kontsevich's obstruction bundle $V_d$ on $\overline M_{0,0}(C,d)$. Its fiber is $H^1(C', f^*(\mathcal O(-1)\oplus \mathcal O(-1)))$. Its top chern class is $\EE_d$. We are now ready to exhibit the connection between the Aspinwall-Morrison calculation and Gromov-Witten invariants.

\begin{proposition}  $\alpha_*\left(e_1^*(h)e_2^*(h)e_3^*(h)\pi^*(\EE_d)\right)=H^{2d+1}.$
\end{proposition} 

{\bf Proof}. Let $E_d$ be the top chern class of the bundle $\rho^*(V_d)$ on $M^1_d$. Recall from part II of the introduction that $H^{2d-2}$ is the top chern class of the Aspinwall-Morrison obstruction bundle $U_d$ on $N^1_d$. It is shown in \cite{[7]} that $\phi_*(E_d)=H^{2d-2}$. On the other hand $\psi^*(E_d)=\pi^*(\EE_d)$. But $\psi$ is birational, hence by the projection formula $\psi_*(\pi^*(\EE_d))=E_d$.

We compute: 
\begin{eqnarray}
& & \alpha_*(\prod_{i}e_i^*(h)\EE_d)=\alpha_*(\prod_{i}e_i^*(h)\psi^*(E_d))=\phi_*(\psi_*(\prod_{i}e_i^*(h))E_d) \nonumber \\ & & =\phi_*(\phi^*(H^3)E_d)=H^3\phi_*(E_d)=H^3H^{2d-2}=H^{2d+1}.
\end{eqnarray} 
The proposition is proven.$\dagger$

\vspace{0.2in}

The last proposition yields:
\begin{equation}
\int_{\overline M_{0,3}(C,d)}\prod_{i=1}^{3}e_i^*(h)\EE_d=\int_{\overline N_d(C)}\alpha_*(\prod_{i=1}^{3}e_i^*(h)\psi^*(E_d))=\int_{\overline N_d(C)}H^{2d+1}=1,
\end{equation}

i.e. the Aspinwall-Morrison calculation is a pushforward of Kontsevich's calculation from $\overline M_{0,3}(C,d)$ to the projective space $\overline N_d(C)$.

E-mail: elezi@math.stanford.edu

Address: Department of Mathematics, Stanford University, Stanford CA, 94305.


\begin{thebibliography}{[Hum]}
\bibitem[1]{[1]} P. S. Aspinwall and D. Morrison, {\it Topological field theory and rational curves}, Comm. Math. Phys., 151(2), (1993), 245-262
\bibitem[2]{[2]} K. Behrend, {\it Gromov-Witten invariants in algebraic geometry}, Invent. Math. {\bf 127} (1997), 601-617.
\bibitem[3]{[13]} K. Behrend and B. Fantechi, {\it The intrinsic normal cone}, Invent. Math., 128 1997, 45-88.
\bibitem[4]{[3]} J. Bryan, S. Katz, N. C. Leung, {\it Multiple covers and the integrality conjecture for rational curves in Calabi-Yau threefolds}, preprint 1999, math.AG/9911056.
\bibitem[5]{[4]} P. Candelas, X. de la Ossa, P. Green and L. Parkes, {\it A pair of Calabi-Yau manifolds as an exactly soluble superconformal theory}, Nuc. Phys. {\bf 539} (1991), 21-74.
\bibitem[6]{[5]} D. Cox and S. Katz, {\it Mirror Symmetry and Algebraic Geometry}, AMS Surveys and Monographs in Mathematics, AMS, Providence RI 1999.
\bibitem[7]{[6]} M. Kontsevich, {\it Enumeration of rational curves via torus actions}, in {\it The moduli space of curves}, (R. Dijkgraaf, C. Faber, and G. van der Geer, eds.), Birkhauser, 1995, 335-168.
\bibitem[8]{[11]} J. Li and G. Tian, {\it Virtual moduli cycles and Gromov-Witten invariants of algebraic varieties},J. AMS 11 1998, 119-174.
\bibitem[9]{[7]} B. Lian, K. Liu, S.-T.Yau, {\it Mirror Principle I}, Asian J. Math. Vol.1, no. 4 (1997), 729-763.
\bibitem[10]{[8]} Yu. I. Mannin, {\it Generating functions in algebraic geometry and sums over trees}, in {\it The moduli of curves}, (R. Dijkgraaf, C. Faber, and G. van der Geer, eds.), Birkhauser, 1995, 401-417.
\bibitem[11]{[9]} R. Pandharipande, {\it Hodge integrals and degenerate contributions}, preprint 1998, math.AG/9811140.
\bibitem[12]{[12]} I. Vainsencher, {\it Enumeration of n-fold tangent hyperplanes to a surface}, J. Algebraic Geom., 4 1995, 503-526.
\bibitem[13]{[10]} C. Voisin, {\it A mathematical proof of a formula of Aspinwall and Morrison}, Compositio Math., 104(2) 1996, 135-151.

\end{thebibliography}
\end{document}